\documentclass[10pt, a4paper]{article}
\usepackage{amsfonts}
\usepackage{mathrsfs}
\usepackage{latexsym}
\usepackage{xy}
\usepackage{amsfonts,amsmath,amssymb,amsthm}
\usepackage{color}

\xyoption{all}

\newcommand{\bcen}{\begin{center}}     \newcommand{\ecen}{\end{center}}
\newcommand{\bay}{\begin{array}}      \newcommand{\eay}{\end{array}}
\newcommand{\beq}{\begin{eqnarray*}}      \newcommand{\eeq}{\end{eqnarray*}}
\def\az{\alpha}
\def\bz{\beta}
\def\dz{\delta}
\def\gz{\gamma}

\def\ot{\otimes}

\def\wt{\widetilde}
\def\mt{\mapsto}

\def\dim{\mathrm{dim}}
\def\dz{\delta}

\def\Hom{\mathrm{Hom}}

\def\mt{\mapsto}

\def\ot{\otimes}

\begin{document}

\newtheorem*{theorem}{Theorem}
\newtheorem{proposition}{Proposition}
\newtheorem{lemma}{Lemma}
\newtheorem{corollary}{Corollary}
\newtheorem*{remark}{Remark}
\newtheorem{example}{Example}
\newtheorem{definition}{Definition}
\newtheorem*{conjecture}{Conjecture}
\newtheorem{question}{Question}

\title{Symmetry criteria for Hochschild extensions}

\author{Yang Han}

\date{\footnotesize KLMM, Academy of Mathematics and Systems Science,
Chinese Academy of Sciences, \\ Beijing 100190, China. \\ School of Mathematical Sciences, University of
Chinese Academy of Sciences, \\ Beijing 100049, China.\\ E-mail: hany@iss.ac.cn}

\maketitle

\begin{abstract}
We give two sufficient and necessary conditions for the Hochschild extension of a finite dimensional algebra by its dual bimodule and a Hochschild 2-cocycle to be a symmetric algebra.
\end{abstract}

\medskip


\medskip



\section{Introduction}

Throughout this paper, $k$ is a fixed field, $\ot :=\ot_k$, and $(-)^* := \Hom_k(-,k)$.
Let $A$ be a finite dimensional $k$-algebra. Then $A^*$ is a natural $A$-bimodule whose $A$-bimodule structure is given by $(cfa)(b):=f(abc)$ for all $f\in A^*$ and $a,b,c\in A$.
The {\it trivial extension} of $A$ is the finite dimensional $k$-algebra $T(A):=A\oplus A^*$ with multiplication $(a,f)(b,g):=(ab,ag+fb)$, which is always a {\it symmetric algebra}, i.e., $T(A)\cong T(A)^*$ as $T(A)$-bimodules, and plays a quite important role in the representation theory of finite dimensional algebras \cite{Hap88}.
Let $\az : A\ot A\to A^*$ be a {\it Hochschild 2-cocycle}, i.e., a $k$-linear map satisfying {\it Hochschild 2-cocycle condition} $a\az(b\ot c)-\az(ab\ot c)+\az(a\ot bc)-\az(a\ot b)c=0$ for all $a,b,c\in A$.
Then the {\it Hochschild extension} $T(A,\az):=A\oplus A^*$ with multiplication
$(a,f)(b,g):=(ab,ag+fb+\az(a\ot b))$ is always a selfinjective algebra \cite[Proposition 1.2]{Yam81},
but unlike trivial extension it is unnecessary to be a symmetric algebra \cite[Proposition 2.8]{OhnTakYam99}.
In \cite[Theorem 2.2]{OhnTakYam99}, Ohnuki-Takeda-Yamagata gave a sufficient condition for a Hochschild extension $T(A,\az)$ to be a symmetric algebra, but the condition is not necessary \cite[Proposition 2.9]{OhnTakYam99}.
In \cite[Theorem 3.1]{Ita19}, Itagaki gave a more general sufficient condition than Ohnuki-Takeda-Yamagata's, but the condition is not necessary either.
Moreover, in \cite[Theorem 2]{HanLiuWan19} we provided a sufficient condition for a Hochschild extension of a finite dimensional dg $k$-algebra to be an $n$-symmetric $A_\infty$-algebra.

In this paper, we will give two sufficient and necessary conditions for a Hochschild extension $T(A,\az)$ to be a symmetric algebra.
Denote by $Z(A)$ the center of a $k$-algebra $A$ and $U(A)$ the group of units (=invertible elements) in $A$. Let\; $\wt{}: \Hom_k(A^{\ot 2},A^*) \to (A^{\ot 3})^*, \az\mt \wt{\az}$, be the canonical adjunction isomorphism where $\wt{\az}$ is given by $\wt{\az}(a\ot b\ot c):=\az(b\ot c)(a)$ for all $a,b,c\in A$. The main result of this paper is the following theorem.

\begin{theorem}
Let $A$ be a finite dimensional $k$-algebra and $\az\in\Hom_k(A\ot A,A^*)$ a Hochschild 2-cocycle. Then the following three statements are equivalent:

{\rm (1)} The Hochschild extension $T(A,\az)$ is a symmetric algebra.

{\rm (2)} There exist $c\in Z(A)\cap U(A)$ and $h\in A^*$ such that
$$\az(a\ot b)(c)-\az(b\ot a)(c)+h(ab-ba)=0$$
for all $a,b\in A$.

{\rm (3)} There exists $c\in Z(A)\cap U(A)$ such that the Hochschild cohomology class $[i_c^*(\wt{\az})]\in H^2(C_\bullet(A)^*)\cong HH_2(A)^*\cong HH^2(A,A^*)$ can be lifted to a cyclic cohomology class in $HC^2(A)\cong HC_2(A)^*$ along the canonical map $I^2: HC^2(A) \to H^2(C_\bullet(A)^*)$ in the Connes periodicity exact sequence.
\end{theorem}

\begin{remark}{\rm
Let $A=kQ/I$ be a bound quiver algebra where $Q$ is a finite quiver and $I$ is an admissible ideal of the path algebra $kQ$ of $Q$ (Ref. \cite{AssSimSkr06}), and $K:=kQ_0$ where $Q_0$ is the vertex set of $Q$. Since $K\cong k^{|Q_0|}$ is separable, the homology of the {\it $K$-relative Hochschild chain complex} of $A$
$$\cdots \xrightarrow{b_3} A \ot_{K^e} A^{\ot_K 2} \xrightarrow{b_2} A \ot_{K^e} A \xrightarrow{b_1} A \ot_{K^e}K \to 0$$
is isomorphic to the Hochschild homology of $A$, and the cohomology of the {\it $K$-relative Hochschild cochain complex of $A$ with coefficients in $A$-bimodule $A^*$}
$$0\to \Hom_{K^e}(K,A^*) \xrightarrow{\dz^0} \Hom_{K^e}(A,A^*) \xrightarrow{\dz^1} \Hom_{K^e}(A^{\ot_K 2},A^*) \xrightarrow{\dz^2} \cdots$$
is isomorphic to the Hochschild cohomology of $A$ with coefficients in $A$-bimodule $A^*$ (Ref. \cite[Section 1]{GerSch86}).
For any Hochschild 2-cocycle $\az\in \Hom_{K^e}(A^{\ot_K 2},A^*)$ and Hochschild extension $T(A,\az)$, we can obtain an almost the same theorem as above in which
the condition (2) with $h=0$ is just Itagaki's sufficient condition in \cite[Theorem 3.1]{Ita19} and
the condition (2) with $h=0$ and $c=1$ is just Ohnuki-Takeda-Yamagata's sufficient condition in \cite[Theorem 2.2]{OhnTakYam99}.
}\end{remark}

\section{Preliminaries}

In this part, we fix some terminologies and notations.

\bigskip

\noindent{\bf Hochschild homology and cyclic homology.} Let $A$ be a (unital associative) $k$-algebra.
The {\it Hochschild chain complex} of $A$ is the $k$-vector space chain complex
$$C_\bullet(A)=(A^{\ot \bullet+1},b_\bullet) :\ \cdots\to A^{\ot n+2} \xrightarrow{b_{n+1}} A^{\ot n+1} \xrightarrow{b_n} A^{\ot n} \to \cdots \to A^{\ot 2} \xrightarrow{b_1} A$$
where the {\it Hochschild boundary operator} $b_n: A^{\ot n+1} \to A^{\ot n}$ is given by
$$\begin{array}{ll}
b_n(a_0\ot a_1\ot\cdots\ot a_n) := &
\sum\limits_{i=0}^{n-1}(-1)^i\ a_0\ot a_1\ot \cdots\ot a_ia_{i+1}\ot\cdots\ot a_n \\ [3mm]
& \quad\quad +(-1)^n\ a_na_0\ot a_1\ot \cdots\ot a_{n-1}.
\end{array}$$
The {\it $n$-th Hochschild homology} $HH_n(A)$ of $A$ is the $n$-th homology of $C_\bullet(A)$.

The {\it cyclic bicomplex} of $A$ is the $k$-vector space bicomplex
$$\xymatrix{
\ar@{.}[d] & \ar@{.}[d] & \ar@{.}[d] & \ar@{.}[d] \\
A^{\ot 4} \ar[d]^-{b_3} & A^{\ot 3} \ar[d]^-{b_2} \ar[l]^-{B_2}
& A^{\ot 2} \ar[d]^-{b_1} \ar[l]^-{B_1} & A \ar[l]^-{B_0} \\
A^{\ot 3} \ar[d]^-{b_2} & A^{\ot 2} \ar[d]^-{b_1} \ar[l]^-{B_1} & A \ar[l]^-{B_0} & \\
A^{\ot 2} \ar[d]^-{b_1} & A \ar[l]^-{B_0} & & \\
A & & &
}$$
where $b_n$ is the Hochschild boundary operator and
$B_n: A^{\ot n+1} \to A^{\ot n+2}$ is the {\it Connes boundary operator} given by
$$\begin{array}{ll}
B_n(a_0\ot\cdots\ot a_n):= & \sum\limits_{i=0}^n(-1)^{ni}\ (1\ot a_i\ot\cdots\ot a_n\ot a_0\ot\cdots\ot a_{i-1} \\ [2mm]
& \quad\quad\quad\quad -a_i\ot 1\ot a_{i+1}\ot\cdots\ot a_n\ot a_0\ot\cdots\ot a_{i-1}).
\end{array}$$
In particular, $B_1(a_0\ot a_1)=1\ot a_0\ot a_1-1\ot a_1\ot a_0-a_0\ot 1\ot a_1+a_1\ot 1\ot a_0$.
The {\it $n$-th cyclic homology} $HC_n(A)$ of $A$ is the $n$-th homology of the total complex of the cyclic bicomplex of $A$.
The {\it Connes periodicity exact sequence} is the long exact sequence
$$\begin{array}{l}
\cdots \xrightarrow{B_{n-1}} HH_n(A) \xrightarrow{I_n} HC_n(A) \xrightarrow{S_n} HC_{n-2}(A) \xrightarrow{B_{n-2}} HH_{n-1}(A) \xrightarrow{I_{n-1}} \\ [2mm]
\cdots \xrightarrow{B_0} HH_1(A) \xrightarrow{I_1} HC_1(A) \xrightarrow{S_1} 0 \xrightarrow{} HH_0(A) \xrightarrow{I_0} HC_0(A) \xrightarrow{} 0 \end{array}$$
where $I_n$ is induced by the injection of Hochschild chain complex in the first column of cyclic bicomplex,
$B_n$ is induced by Connes boundary operator, and $S_n$ is the {\it periodicity map} induced by shifting (Ref. \cite[2.2.1]{Lod98}).

\bigskip

\noindent{\bf Hochschild cohomology and cyclic cohomology.}
Let $A$ be a $k$-algebra and $M$ an $A$-bimodule. The {\it Hochschild cochain complex of $A$ with coefficients in $M$} is the $k$-vector space cochain complex $C^\bullet(A,M)=(\Hom_k(A^{\ot \bullet},M),\dz^\bullet)$:
$$0\to M \xrightarrow{\dz^0} \Hom_k(A,M) \xrightarrow{\dz^1} \Hom_k(A^{\ot 2},M) \xrightarrow{\dz^2} \Hom_k(A^{\ot 3},M) \to \cdots$$
where the {\it Hochschild coboundary operators} $\dz^0$ is given by $\dz^0(m)(a):=am-ma$ and $\dz^n: \Hom_k(A^{\ot n},M) \to \Hom_k(A^{\ot n+1},M)$ is given by
$\dz^n(\az)(a_1\ot\cdots\ot a_{n+1}) := a_1\az(a_2\ot\cdots\ot a_{n+1})+
\sum\limits_{i=1}^n(-1)^i\ \az(a_1\ot\cdots\ot a_ia_{i+1}\ot\cdots\ot a_{n+1})+(-1)^{n+1}\ \az(a_1\ot\cdots\ot a_n)a_{n+1}$ for all $n\ge 1$.
In particular,
$\dz^2(\az)(a_1\ot a_2\ot a_3)=a_1\az(a_2\ot a_3)-\az(a_1a_2\ot a_3)+\az(a_1\ot a_2a_3)-\az(a_1\ot a_2)a_3$ for all $a_1,a_2,a_3\in A$.
The {\it $n$-th Hochschild cohomology $HH^n(A,M)$ of $A$ with coefficients in $M$} is
the $n$-th cohomology of the Hochschild cochain complex $C^\bullet(A,M)$ of $A$ with coefficients in $M$. The {\it $n$-th Hochschild cohomology} $HH^n(A)$ of $A$ is just $HH^n(A,A)$.

Applying the $k$-dual functor $(-)^*$ to the Hochschild chain complex $C_\bullet(A)=(A^{\ot \bullet+1},b_\bullet)$ of $A$,
we obtain the cochain complex $C_\bullet(A)^*=((A^{\ot \bullet+1})^*,b_\bullet^*)$.
According to the canonical adjunction isomorphism $(A^{\ot n+1})^* \cong \Hom_k(A^{\ot n},A^*)$,
the Hochschild cochain complex $C^\bullet(A,A^*)=(\Hom_k(A^{\ot \bullet},A^*),\dz^\bullet)$ of $A$ with coefficients in $A^*$
is isomorphic to the cochain complex $C_\bullet(A)^*=((A^{\ot \bullet+1})^*,b_\bullet^*)$.
In particular, $H^n(C_\bullet(A)^*)\cong HH_n(A)^*\cong HH^n(A,A^*)$ for all $n\in \mathbb{N}_0$.

Applying the $k$-dual functor $(-)^*$ to the cyclic bicomplex of $A$,
we obtain the {\it cyclic cochain bicomplex} of $A$
$$\xymatrix{
& & & \\
(A^{\ot 4})^* \ar[r]^-{B_2^*} \ar@{.}[u] & (A^{\ot 3})^* \ar[r]^-{B_1^*} \ar@{.}[u]
& (A^{\ot 2})^* \ar[r]^-{B_0^*} \ar@{.}[u] & A^* \ar@{.}[u] \\
(A^{\ot 3})^* \ar[u]^-{b_3^*} \ar[r]^-{B_1^*} & (A^{\ot 2})^* \ar[u]^-{b_2^*} \ar[r]^-{B_0^*} & A^* \ar[u]^-{b_1^*} & \\
(A^{\ot 2})^* \ar[u]^-{b_2^*} \ar[r]^-{B_0^*} & A^* \ar[u]^-{b_1^*} & & \\
A^* \ar[u]^-{b_1^*} & & & \\
}$$
The {\it $n$-th cyclic cohomology} $HC^n(A)$ of $A$ is the $n$-th cohomology of the total complex of the cyclic cochain bicomplex of $A$ (Ref. \cite[2.4]{Lod98}). Clearly, $HC^n(A) \cong HC_n(A)^*$ for all $n \in \mathbb{N}_0$. From the cyclic cochain bicomplex of $A$ we can obtain the {\it Connes periodicity exact sequence} (cohomology form)
$$\begin{array}{l}
0 \xrightarrow{} HC^0(A) \xrightarrow{I^0} H^0(C_\bullet(A)^*) \xrightarrow{} 0 \xrightarrow{S^1} HC^1(A)
\xrightarrow{I^1} H^1(C_\bullet(A)^*) \xrightarrow{B^0} \cdots \\ [2mm]
\xrightarrow{I^{n-1}} H^{n-1}(C_\bullet(A)^*) \xrightarrow{B^{n-2}} HC^{n-2}(A) \xrightarrow{S^n} HC^n(A) \xrightarrow{I^n} H^n(C_\bullet(A)^*) \xrightarrow{} \cdots
\end{array}$$
where $I^n$ is induced by the projection of the cyclic cochain bicomplex to the first column,
$B^n$ is induced by the dual $B_n^*$ of the Connes boundary operator $B_n$, and $S^n$ is the {\it periodicity map} induced by shifting. In particular, $I^2([(\bz,\gz)])=[\bz]$ where $(\bz,\gz)\in (A^{\ot 3})^*\oplus A^*$ is any cyclic 2-cocycle and $\bz\in (A^{\ot 3})^*$ is a Hochschild 2-cocycle.

\bigskip

\noindent{\bf Hochschild extensions.} Let $A$ be a $k$-algebra, $M$ an $A$-bimodule,
and $\az : A \ot A \to M$ a {\it Hochschild 2-cocycle},
i.e., a $k$-linear map satisfying the {\it Hochschild 2-cocycle condition} $a\az(b\ot c)-\az(ab\ot c)+\az(a\ot bc)-\az(a\ot b)c=0$ for all $a,b,c\in A$.
The {\it Hochschild extension} $T(A,M,\az)$ of $A$ by $M$ and $\az$ is
the $k$-algebra $A\oplus M$ with multiplication $(a,m)(a',m'):=(aa',am'+ma'+\az(a\ot a'))$ (Ref. \cite[1.5.3]{Lod98}).
In the case of $M=A^*$, we denote by $T(A,\az)$ the Hochschild extension $T(A,A^*,\az)$ for short.

\bigskip

\noindent {\bf Contraction operator.} Let $A$ be a $k$-algebra. For any Hochschild $m$-cochain ${\az}\in C^m(A):=C^m(A,A)$, the {\it contraction operator} $i_{\az}$ on the Hochschild chain complex $C_\bullet(A)$ is the graded morphism $i_{\az}: C_\bullet(A) \to C_\bullet(A)$ of degree $-m$ given by
$$i_{\az}(a_0\ot a_1\ot \cdots\ot a_n) := a_0{\az}(a_1\ot\cdots\ot a_m) \ot (a_{m+1}\ot\cdots\ot a_n).$$
Note that $bi_{\az}-i_{\az}b= i_{\dz\az}$ where $b$ is the Hochschild boundary operator and $\dz$ is the Hochschild coboundary operator.
Consequently, if $\az$ is a Hochschild cocycle then the contraction operator $i_{\az}$ on $C_\bullet(A)$ induces an operator $i_{\az}$ on $HH_\bullet(A)$.
In particular, for any $c\in Z(A)=HH^0(A)\subseteq C^0(A)=A$, $i_c$ is an operator on $C_\bullet(A)$ which is given by $i_c(a_0\ot a_1\ot \cdots \ot a_n):=a_0c\ot a_1\ot\cdots\ot a_n$.
Moreover, $b_ni_c=i_cb_n$ for all $n\ge 1$.
For more properties of contraction operators in noncommutative differential calculus theory, we refer to \cite{TamTsy05,HanLiuWan18}.

\section{The proof of Theorem}

Now we start the proof of Theorem.
Write $T(A,\az)$ as $T$ for short.
Any element in $T^*=\Hom_k(A\oplus A^*,k) \cong A^*\oplus A^{**} \cong A^*\oplus A$ can be written as $(f,a)$ where $f\in A^*, a\in A$, and the $k$-linear form $(f,a): T \to k$ is given by $(f,a)((b,g)):=f(b)+g(a)$ for all $(b,g)\in T $.

\bigskip

\noindent{\bf (1)$\Rightarrow$(2):}
Assume that $T$ is a symmetric algebra. Then there exists a $T$-bimodule isomorphism $\phi: T \to T^*$.
Let $\phi(1_{T})=(h,c)\in T^*$ where $h\in A^*$ and $c\in A$.

\medskip

{\bf Claim 1.} {\it For all $(a,f), (b,g)\in T$, we have
$$\az(a\ot b)(c)-\az(b\ot a)(c)+h(ab-ba)=0,\ f(bc)=f(cb),\ g(ac)=g(ca).$$}
Indeed, since $\phi$ is a $T$-bimodule morphism, for any $(a,f)\in T$, we have
$$\phi(a,f)=(a,f)\phi(1_{T})=\phi(1_{T})(a,f)$$
which maps any $(b,g)\in T$ to
$$\phi(a,f)((b,g))=((a,f)\phi(1_{T}))((b,g))=(\phi(1_{T})(a,f))((b,g)).$$
By the canonical $T$-bimodule structure of $T^*$, we get
$$\phi(a,f)((b,g))=\phi(1_{T})((b,g)(a,f))=\phi(1_{T})((a,f)(b,g)).$$
According to the multiplication of $T$, we have
$$\phi(a,f)((b,g))=(h,c)((ba,bf+ga+\az(b\ot a)))=(h,c)((ab,ag+fb+\az(a\ot b))).$$
From the definition of $(h,c)\in T^*$, we obtain
$$\begin{array}{ll}
\phi(a,f)((b,g)) & =h(ba)+(bf)(c)+(ga)(c)+\az(b\ot a)(c) \\
& =h(ab)+(ag)(c)+(fb)(c)+\az(a\ot b)(c).
\end{array}$$
By the canonical $T$-bimodule structure of $T^*$, we get
$$\begin{array}{ll}
\phi(a,f)((b,g)) & =h(ba)+f(cb)+g(ac)+\az(b\ot a)(c) \\
& =h(ab)+g(ca)+f(bc)+\az(a\ot b)(c).
\end{array}$$
Taking $a=0$, we obtain $f(bc)=f(cb)$ for all $b\in A$ and $f\in A^*$.
Taking $b=0$, we get $g(ac)=g(ca)$ for all $a\in A$ and $g\in A^*$.
Furthermore, for all $(a,f), (b,g)\in T$, we have
$\az(a\ot b)(c)-\az(b\ot a)(c)+h(ab-ba)=0.$

\medskip

{\bf Claim 2.} $c\in Z(A)$.

By Claim 1, we have $f(bc-cb)=0$ for all $b\in A$ and $f\in A^*$.
Thus $bc-cb=0$ for all $b\in A$. Hence $c\in Z(A)$.

\medskip

{\bf Claim 3.} $c\in U(A)$.

Assume on the contrary $c\notin U(A)$.
It follows from Claim 2 that the two-sided ideal $AcA$ of $A$ generated by $c$ is equal to both $Ac$ and $cA$. However, it is not equal to $A$ due to $c\notin U(A)$.
Thus the complement space $(Ac)^\perp$ of $Ac$ in $A$ is not zero.
Since $A=Ac\oplus (Ac)^\perp$, we have $A^*=(Ac)^*\oplus ((Ac)^\perp)^*$.
Thus there exists $f\in ((Ac)^\perp)^* \subseteq A^*$ such that $f\ne 0$ and $f(Ac)=0$.
For any $(b,g)\in T$, as in the proof of Claim 1, we have $\phi(0,f)((b,g))=h(0b)+g(c0)+f(bc)+\az(0\ot b)(c)=f(bc)\in f(Ac)=0$. So $\phi(0,f)=0$.
Since $\phi$ is a $T$-bimodule isomorphism, we have $(0,f)=(0,0)$, and further $f=0$.
It is a contradiction.

\medskip

Finally, it follows from Claim 1, 2 and 3 that there exist $c\in Z(A)\cap U(A)$ and $h\in A^*$
such that $\az(a\ot b)(c)-\az(b\ot a)(c)+h(ab-ba)=0$ for all $a,b\in A$.

\bigskip

\noindent{\bf (2)$\Rightarrow$(3):} Assume that there exists $c\in Z(A)\cap U(A)$ such that
$\az(a\ot b)(c)-\az(b\ot a)(c) + h(ab-ba)=0$ for all $a,b\in A$.
From Hochschild 2-cocycle condition, we obtain $\az(1\ot ab)=\az(1\ot a)b$ for all $a,b\in A$.
Due to $c\in Z(A)$, we have
$$\begin{array}{ll}
& (B_1^*(i_c^*(\wt{\az})))(a\ot b) \\
= & (\wt{\az}i_cB_1)(a\ot b) \\
= & (\wt{\az}i_c)(1\ot a\ot b-1\ot b\ot a-a\ot 1\ot b+b\ot 1\ot a) \\
= & \wt{\az}(c\ot a\ot b-c\ot b\ot a-ac\ot 1\ot b+bc\ot 1\ot a) \\
= & \wt{\az}(c\ot a\ot b)-\wt{\az}(c\ot b\ot a)-\wt{\az}(ac\ot 1\ot b)+\wt{\az}(bc\ot 1\ot a) \\
= & \az(a\ot b)(c)-\az(b\ot a)(c)-\az(1\ot b)(ac)+\az(1\ot a)(bc) \\
= & -h(ab-ba)+(\az(1\ot 1)c)(ab-ba) \\
= & -(h-\az(1\ot 1)c)(ab-ba) \\
= &-((h-\az(1\ot 1)c)b_1)(a\ot b) \\
= & -b_1^*(h-\az(1\ot 1)c)(a\ot b).
\end{array}$$
Let $h'=h-\az(1\ot 1)c \in A^*$. Then $B_1^*(i_c^*(\wt{\az}))+b_1^*(h')=0$.
Thus the Hochschild cohomology class $[i_c^*(\wt{\az})]\in H^2(C_\bullet(A)^*)$ can be lifted to
the cyclic cohomology class $[(i_c^*(\wt{\az}),h')]\in HC^2(A)$
along the canonical map $I^2:HC^2(A)\to H^2(C_\bullet(A)^*)$.

\bigskip

\noindent{\bf (3)$\Rightarrow$(2):}
Assume that there exists $c\in Z(A)\cap U(A)$ such that the Hochschild cohomology class
$[i_c^*(\wt{\az})] \in H^2(C_\bullet(A)^*)$ can be lifted to a cyclic cohomology class
$[(\az',h')]\in HC^2(A)$ along the canonical map $I^2:HC^2(A)\to H^2(C_\bullet(A)^*)$ where $\az'\in(A^{\ot 3})^*$ and $h'\in A^*$, i.e., $I^2([(\az',h')])=[\az']=[i_c^*(\wt{\az})]$.
Then $b_3^*(\az')=0$, $B_1^*(\az')+b_1^*(h')=0$, i.e., $\az'B_1+h'b_1=0$, and $[\az']=[i_c^*(\wt{\az})] \in H^2(C_\bullet(A)^*)$.
Thus there is $\bz \in (A^{\ot 2})^*$ such that $i_c^*(\wt{\az})=\az'+b_2^*(\bz)$, i.e., $\wt{\az}i_c=\az'+\bz b_2$.
So $\wt{\az}i_cB_1=\az'B_1+\bz b_2B_1=-h'b_1-\bz B_0b_1=-(h'+\bz B_0)b_1$.
As in the proof of (2)$\Rightarrow$(3), we have
$\az(a\ot b)(c)-\az(b\ot a)(c)+(h'+\bz B_0 +\az(1\ot 1)c)(ab-ba)=0$ for all $a,b\in A$.
Let $h=h'+\bz B_0 +\az(1\ot 1)c\in A^*$.
Then $\az(a\ot b)(c)-\az(b\ot a)(c)+h(ab-ba)=0$ for all $a,b\in A$.

\bigskip

\noindent{\bf (2)$\Rightarrow$(1):}
Assume that there exists $c\in Z(A)\cap U(A)$ and $h\in A^*$ such that
$\az(a\ot b)(c)-\az(b\ot a)(c)+h(ab-ba)=0$ for all $a,b\in A$.
Define a map $\phi : T\to T^*$ by $\phi(a,f)((b,g)) := (h,c)((a,f)(b,g))$.
Note that
$(h,c)((a,f)(b,g)) =h(ab)+g(ca)+f(bc)+\az(a,b)(c)$
and
$(h,c)((b,g)(a,f)) =h(ba)+f(cb)+g(ac)+\az(b,a)(c).$
Due to $c\in Z(A)$, we get $f(bc)=f(cb)$ and $g(ac)=g(ca)$.
Furthermore, for all $(a,f),(b,g)\in T$, we have
$$\phi(a,f)((b,g)) = (h,c)((a,f)(b,g))=(h,c)((b,g)(a,f)),$$
or equivalently,
$$\phi(a,f)=(h,c)(a,f) =(a,f)(h,c).$$

Now we prove that $\phi$ is a $T$-bimodule morphism.
Indeed, for all $(a',f'),\linebreak (a,f),(a'',f'')\in T$, we have
$$\begin{array}{ll}
\phi((a',f')(a,f)(a'',f'')) & =(h,c)(a',f')(a,f)(a'',f'') \\
& =(a',f')(h,c)(a,f)(a'',f'') \\
& =(a',f')\phi(a,f)(a'',f'').
\end{array}$$

Next we show that $\phi$ is a $T$-bimodule isomorphism.
Since $\dim_k T=\dim_k T^*<\infty$, it suffices to prove that $\phi$ is injective.
For this, suppose $\phi(a,f)=0$. We need to prove $(a,f)=0$.
Indeed, on one hand, for any $g\in A^*$, we have $0 = \phi(a,f)((0,g))= h(a0)+g(ca)+f(0c)+\az(a\ot 0)(c)=g(ca)$. Thus $ca=0$.
Due to $c\in U(A)$, we get $a=0$.
On the other hand, for any $b\in A$, we have $0 = \phi(0,f)((b,0)) = h(0b)+0(c0)+f(bc)+\az(0\ot b)(c)=f(bc)$.
Due to $c\in U(A)$, we obtain $f(A)=0$. Thus $f=0$.

\bigskip

Now we finish the proof of Theorem.

\bigskip

\noindent {\footnotesize {\bf ACKNOWLEDGEMENT.} The author is sponsored by Project 11971460 NSFC.}


\footnotesize

\begin{thebibliography}{99}

\bibitem{AssSimSkr06} I. Assem, D. Simson and A. Skowro\'{n}ski, Elements of the representation theory of associative algebras, Vol. 1, Techniques of representation theory, London Mathematical Society Student Texts 65, Cambridge University Press, Cambridge, 2006.

\bibitem{GerSch86} M. Gerstenhaber and S.D. Schack, Relative Hochschild cohomology, rigid algebras,
and the Bockstein, J. Pure Appl. Algebra 43 (1986), 53--74.

\bibitem{HanLiuWan18} Y. Han, X. Liu and K. Wang, Hochschild (co)homologies of dg $K$-rings and their Koszul duals, arXiv:1810.05969v1 [math.KT].

\bibitem{HanLiuWan19} Y. Han, X. Liu and K. Wang, Exact Hochschild extensions and deformed Calabi-Yau completions, arXiv:1909.02200v1 [math.RA].

\bibitem{Hap88} D. Happel, Triangulated categories in the representation theory of finite dimensional algebras, London Math. Soc. Lecture Notes Series 119, Cambridge University Press, 1988.

\bibitem{Ita19} T. Itagaki, Symmetric Hochschild extension algebras and normalized 2-cocycles, Arch. Math. 112 (2019), no. 3, 249--259.

\bibitem{Lod98} J.-L. Loday, Cyclic homology, The second edition,
Grundlehren der mathematischen Wissenschaften 301, Springer-Verlag, Berlin Heidelberg, 1998.

\bibitem{OhnTakYam99} Y. Ohnuki, K. Takeda and K. Yamagata, Symmetric Hochschild extension algebras, Colloq. Math. 80 (1999), 155--174.

\bibitem{TamTsy05} D. Tamarkin and B. Tsygan, The ring of differential operators on forms in noncommutative calculus, In: Graphs and patterns in mathematics and theoretical physics, Proc. Sympos. Pure Math. 73, Amer. Math. Soc., Providence, RI 2005, 105--131.

\bibitem{Yam81} K. Yamagata, Extensions over hereditary artinian rings with self-dualities, I,
J. Algebra 73 (1981), 386--433.

\end{thebibliography}
\end{document}